\documentclass[12pt]{amsart}
\usepackage{amsfonts}
\usepackage{amssymb}
\usepackage{graphicx}

               \def\z{\zeta}

\def\D{{\mathbb D}}  
\def\C{{\mathbb C}}  
\def\R{{\mathbb R}}

\def\({\left(}       \def\){\right)}

\newtheorem{lem}{\sc Lemma}
\newtheorem{thm}{\sc Theorem}

\newtheorem{other}{\sc Theorem}              
\newenvironment{pf}{\noindent{\textit{Proof. }}}{$\Box$ }

\title[Bounds for the growth and distortion]{Sharp bounds for the growth and distortion of the analytic part of convex harmonic functions}

\author[María J. Martín]{Mar\'{\i}a J. Mart\'{\i}n}
\address{Departamento de An\'alisis Matemático \& IMAULL, Universidad de La Laguna. Av.  Astrofísico Francisco Sánchez, s/n. Facultad de Matemáticas. 38200, La Laguna, Tenerife, Spain.} \email{maria.martin@ull.es}

\date{\today}

\begin{document}


\begin{abstract} We obtain the sharp upper and lower bounds for the growth and distortion of the analytic parts $h$ of orientation-preserv\-ing harmonic mappings  $f=h+\overline g$ (normalized in the standard way) that map the unit disk onto a convex domain. \end{abstract}

\maketitle
\section{Introduction}
A locally univalent,  orientation-preserving, complex-valued harmonic function $f=h+\overline g$ in the unit disk $\D$, normalized by the conditions $h(0)=g(0)=1-h'(0)=0$ is said to belong to the class $\mathcal K_H$ if $f(\D)$ is a convex domain. If, in addition, $g'(0)=0$, the function $f$ belongs to the class $\mathcal K_H^0$. Notice that both functions $h$ and $g$ in the decomposition $f=h+\overline g$ are holomorphic in the unit disk whenever $f$ is a harmonic function. 

\par\smallskip
The obvious fact that any analytic function is harmonic, shows that the class $\mathcal K$, which consists of those analytic functions $\varphi$ in $\D$ such that $\varphi(0)=1-\varphi'(0)=0$, is contained in $\mathcal K_H^0$. 

\par\smallskip
It is well known (see \cite[Thm. 2.2.8] {GK}, for instance) that the growth and distortion of a given $\varphi\in \mathcal K$ are controlled, respectively, by the estimates
\begin{equation}\label{eq-growth-K}
\frac{|z|}{1+|z|}\leq |\varphi(z)|\leq \frac{|z|}{1-|z|}\,,\quad z\in\D\,,
\end{equation}
and 
\begin{equation}\label{eq-dist-K}
\frac{1}{(1+|z|)^2}\leq |\varphi'(z)|\leq \frac{1}{(1-|z|)^2}\,,\quad z\in\D\,.
\end{equation}
 In particular, if we let $|z|\to 1$ in the lower bound in \eqref{eq-growth-K}, we conclude that the Euclidean disk of center $0$ and radius $1/2$ is contained in the range of any function $\varphi\in \mathcal K$. The value $1/2$ is sharp, as will be justified below. 
 
\par\smallskip
Both \eqref{eq-growth-K} and \eqref{eq-dist-K} are sharp: equality in any of the inequalities mentioned is obtained, for a given $z\in\D\setminus\{0\}$ by considering the function $\ell$, defined by
\begin{equation}\label{eq-ell}
\ell(z)=\frac{z}{1-z}\,,\quad z\in\D\,,
\end{equation}
or some of its rotations, $\ell_\lambda$, where, for $|\lambda|=1$ and $z$ in the unit disk,  $\ell_\lambda(z)=\overline\lambda \ell(\lambda z)$\,. 

\par\smallskip
The function $\ell$ maps the unit disk onto the half-plane $\mathbb H=\{z\in\C\colon {\rm Re\,}z >-1/2\}$. Therefore, no disk centered at the origin and with radius $r>1/2$ is contained in the range of $\ell$.
\par\smallskip
The class $\mathcal K_H^0$ is much wider that its analytic counterpart $\mathcal K$ and there is a number of difficult problems that remain unresolved in relation with functions in this class. For instance, no sharp bounds for the growth of a given $f\in \mathcal K_H^0$ are known. More concretely, we mention that it is not difficult to check (see \cite[p. 100]{Dur-harm}) that if $f\in \mathcal K_H^0$ and $z\in\D$, then
\begin{equation}\label{eq-growth-Ha}
\frac{|z|}{(1+|z|)^2}\leq |f(z)|\leq \frac{|z|}{(1-|z|)^2}\,.
\end{equation}

Nevertheless, as pointed out in \cite[p. 100]{Dur-harm}, the lower bound in \eqref{eq-growth-Ha} is far from being sharp. The  upper bound is known to be correct with respect to the order of growth since
\[
\limsup_{r\to 1}\ (1-r)^2 \sup_{|z|=r} |L(z)|\geq \frac{1}{2}\,,
\]
where $L=H+\overline G \in \mathcal K_H^0$ is the \emph{half-plane harmonic mapping} defined by
\begin{equation}\label{eq-L}
L(z)=\frac{1}{2}\left[\left(\ell(z)+k(z)\right)+ \left(\overline{\ell(z)-k(z)}\right)\right]\,,\quad z\in\D\,.
\end{equation}
In \eqref{eq-L},  $\ell$ is given by \eqref{eq-ell} and $k$ is the \emph{Koebe function}
\begin{equation}\label{eq-koebe}
k(z)=\frac{z}{(1-z)^2}\,,\quad z\in\D\,.
\end{equation}
Notice that the \emph{analytic part} $H$ of $L$ equals to
\begin{equation}\label{eq-H}
H(z)=\frac{2z-z^2}{2(1-z)^2}\,,\quad z\in\D\,.
\end{equation}

\par\smallskip
We again refer the reader to  \cite[p. 100]{Dur-harm}, where the author states that it seems likely that the sharp bounds for the growth of functions in the class $\mathcal K_H^0$ are attained by the function $L$. This is a problem that is still open and will not be solved in this paper. However, as a perhaps first step towards a solution, we obtain the sharp bounds for the growth and distortion of the analytic part $h$ of functions $f=h+\overline g$ in the class $\mathcal K_H^0$, a result which seems to be interesting by itself and whose proof seems to have eluded the experts in the area so far.  Here is the specific statement.

\begin{thm}\label{thm-main}
Let $f=h+\overline g \in \mathcal K_H^0$. Then, the growth and distortion of the analytic part $h$ of $f$ are subject, respectively, to the bounds
\begin{equation}\label{eq-growth-H}
 \frac{2|z|+|z|^2}{2(1+|z|)^2}\leq |h(z)|\leq \frac{2|z|-|z|^2}{2(1-|z|)^2}\,,\quad z\in \D\,,
\end{equation}
and 
\begin{equation}\label{eq-dist-H}
\frac{1}{(1+|z|)^3}\leq |h'(z)|\leq \frac{1}{(1-|z|)^3}\,,\quad z\in \D\,.
\end{equation}
The estimates are sharp: equality holds in any of the inequalities if and only if $f$ equals the harmonic half-plane mapping \eqref{eq-L} or some of its rotations $L_\lambda$ given by $L_\lambda(z)=\overline\lambda L(\lambda z)$, where $|\lambda|=1$\,.  That is, if and only if $h$ equals the function $H$ in \eqref{eq-H} or some of its rotations.
\end{thm}

As a consequence of our theorem, if we let $|z|\to 1$ in the lower bound in \eqref{eq-growth-H}, we conclude that the Euclidean disk of center $0$ and radius $3/8$ is contained in the range of every function $h$ which is the analytic part of a function  $f\in \mathcal K_H^0$. The value $3/8$ is sharp, since the function $H$ given by \eqref{eq-H} satisfies that no Euclidean disk centered at the origin with radius $r>3/8$ is contained in $H(\D)$.

\section{Background}\label{sec-back}
In this section we review some results that will be used in the proof of Theorem~\ref{thm-main}.

\subsection{Growth and distortion of functions in the class $\mathcal S$}\label{sec-S} The class $\mathcal S$ consists on univalent (\emph{i.e.} one-to-one) analytic functions $\varphi$ in the unit disk normalized by the conditions $\varphi(0)=1-\varphi'(0)=0$. The following result is well known. We refer the reader to \cite[Thms. 2.5, 2.6]{Dur-univ}, \cite[Thm. 2.2.7]{GK}, or \cite[Thm. 1.6]{Pom-univ} for a proof. 
\begin{other}\label{thm-S}
Let $\varphi\in\mathcal S$ and let $z\in\D$. Then, 
\begin{equation*}\label{eq-growth-S}
\frac{|z|}{(1+|z|)^2}\leq |\varphi(z)|\leq \frac{|z|}{(1-|z|)^2}
\end{equation*}
and 
\begin{equation*}\label{eq-dist-S}
\frac{1-|z|}{(1+|z|)^3}\leq |\varphi'(z)|\leq \frac{1+|z|}{(1-|z|)^3}\,.
\end{equation*}
Equality in any of the inequalities mentioned is obtained, for a given $z\in\D\setminus\{0\}$, by the Koebe function $k$ in \eqref{eq-koebe} or some of its rotations\,. 
\end{other}

Another relevant theorem that must be mentioned is due to Bieberbach \cite{B1, B2} (see also \cite[Thm. 2.2]{Dur-univ}).
\begin{other}\label{thm-Bie}
Let $\varphi\in\mathcal S$. Then, $|\varphi''(0)|\leq 4$. Moreover, $|\varphi''(0)|= 4$ if and only if $\varphi$ equals a rotation of the Koebe function $k$.
\end{other}

\subsection{Functions in the Carath\'eodory class} Let $p$ be an analytic function in the unit disk with $p(0)=1$. If ${\rm Re\,} p(z)\geq 0$ for $z\in\D$, we say that $p$ belongs to the \emph{Carath\'eodory class}, denoted by $\mathcal P$.

\par\smallskip
Any $p\in\mathcal P$ can be represented as
\[
p=\frac{1+\delta}{1-\delta}
\]
for some analytic function $\delta$ in the unit disk such that $\delta(0)=0$ and $\delta(\D)\subset\D$ (see \cite[p. 28]{GK}). As a direct consequence of the Schwarz lemma, we have that if $p\in\mathcal P$,  
\begin{equation}\label{eq-growth-P}
\frac{1-|z|}{1+|z|}\leq |p(z)|\leq \frac{1+|z|}{1-|z|}\,,\quad z\in\D\,.
\end{equation}
The inequalities in \eqref{eq-growth-P} are sharp, as rotations of the functions $p(z)=(1+z)/(1-z)$, which belong to $\mathcal P$, show. 

\par\smallskip
The following lemma will be used in Section~\ref{sec-main} to prove Theorem~\ref{thm-main}. The argument used in the proof can be found in \cite{R}.
\begin{lem}\label{lem-q}
Let $q$ be an analytic function in the unit disk with $q(0)=1$. Suppose that there exists $\alpha\in\R$ such that ${\rm Re\,} \left(e^{i\alpha} q(z)\right)\geq 0$ for $z\in\D$. Then, there exists a function $\delta$, analytic in the unit disk, with $\delta(0)=0$ and $\delta(\D)\subset \D$ such that
\begin{equation}\label{eq-q}
q(z)=\frac{1+e^{-2i\alpha}\delta(z)}{1-\delta(z)}\,,\quad z\in\D\,.
\end{equation}
\end{lem}
\begin{pf}
Since ${\rm Re\,} \left(e^{i\alpha} q(z)\right)>0$ and $q(0)=1$, we have that $\cos\alpha>0$. Let 
\[
Q(z)=\frac{e^{i\alpha}q(z)-i\sin\alpha}{\cos\alpha}\,,\quad z\in\D\,.
\]
Then, $Q$ is analytic in the unit disk, $Q(0)=1$, and ${\rm Re\,} Q(z)>0$ for all $z$ in the unit disk. That is, $Q\in\mathcal P$. Therefore, there exists a function $\delta$ as in the statement of the theorem such that  $Q=(1+\delta)/(1-\delta)$. This gives \eqref{eq-q} and ends the proof of the lemma. 
\end{pf}

\subsection{Functions convex in some direction}\label{sec-convexinsomedirection}
A domain $\Omega$ in the complex plane is said to be \emph{convex in the direction} $t\in[0, \pi)$ if the intersection of $\Omega$ with any line parallel to the line through $0$ and $e^{it}$ is an interval or empty. Note that a convex domain is a domain convex in every direction. The orientation-preserving complex-valued function $f$ in the unit disk $\D$ is said to be convex in the $t$ direction if $f(\D)$ is convex in the direction~$t$. Since every domain convex in some direction is a close-to-convex domain, any locally univalent analytic function convex in some direction is univalent (see, for instance, \cite{Pom-CtC} and the references therein).

\par\smallskip

\subsection{Some fundamental  results on convex harmonic mappings}

The following beautiful theorem was proved in \cite{C-SS}. It shows a relation between functions in  $\mathcal K_H^0$ and the classes of functions convex in some direction defined in Section~\ref{sec-convexinsomedirection} and the class $\mathcal S$ introduced in Section~\ref{sec-S}.

\begin{other}\label{thm-CSS1}
Let $f=h+\overline g$ be a locally univalent,  orientation-preserv\-ing harmonic function with $h(0)=g(0)=1-h'(0)=g'(0)=0$. Then, $f\in  \mathcal K_H^0$  if and only if for all $\theta\in[0,\pi)$ the analytic function $h-e^{2i\theta}g$ is convex in the direction $\theta$.   In particular, if $f=h+\overline g \in\mathcal K_H^0$, then for all $|\varepsilon|=1$ the analytic functions $h-\varepsilon g$ belong to the class $\mathcal S$.
\end{other}

\par\smallskip
Another  result from \cite{C-SS} is the next one (\emph{cf.} \cite[Lem. 5.11]{C-SS}).

\begin{other}\label{thm-CSS2}
Let $f=h+\overline g\in \mathcal K_H^0$.  Then, there exist $\alpha, \beta\in\R$ such that for all $z\in\D$,
\begin{equation*}\label{eq-lema}
{\rm Re\,}\left( e^{i(\alpha+\beta)} \left(h'(z)+e^{-2i\alpha}g'(z)\right)\left(1-e^{-2i\beta}z^2\right)\right)\geq 0\,.
\end{equation*}
\end{other}

Regarding the Taylor coefficients of the analytic part $h$ and the anti-analytic part $g$ of functions $f=h+\overline g$ in $\mathcal K_H^0$, we recall the result in the following theorem, whose proof can be found in \cite{C-SS} or in \cite[p. 50]{Dur-harm}. 

\begin{other}\label{thm-CSS3}
Let $f=h+\overline g\in \mathcal K_H^0$.  Assume that the analytic functions $h$ and $g$ have Taylor series expansions
\[
h(z)=z+\sum_{n=2}^\infty a_nz^n \quad \text{and}\quad g(z)=\sum_{n=2}^\infty b_nz^n\,,\quad z\in\D\,, 
\]
respectively. Then,  for $n=2, 3, \ldots$,
\[
|a_n|\leq \frac{n+1}{2}\quad \text{and}\quad |b_n|\leq \frac{n-1}{2}\,.
\]
Equality holds for the function $L$ defined in \eqref{eq-L}.
\end{other}

\par\smallskip
The \emph{dilatation} of a locally-univalent, orientation-preserving harmonic function $f=h+\overline g$ in $\D$ is the analytic function in the unit disk defined by $\omega=g'/h'$. The fact that this function is analytic in the unit disk follows from the fact that de Jacobian $J_f=|h'|^2-|g'|^2$ of any locally univalent orientation-preserving harmonic mapping $f=h+\overline g$ is positive in $\D$ (see \cite{Lewy}). Indeed this condition implies that not only $\omega$ is holomorphic in the unit disk but also that  $\omega(\D)\subset\D$.  It is clear that if $f\in \mathcal K_H^0$ has dilatation $\omega$, then $\omega(0)=0$. 

\par\smallskip
To finish this section, we mention a transformation that preserves the class $\mathcal K_H^0$. It is related to rotations, which have been already mentioned. Concretely, as pointed out above, given $f=h+\overline g \in\mathcal K_H^0$ and given $|\lambda|=1$, the \emph{rotation} $f_\lambda$ of $f$ is defined, for $z$ in the unit disk, by $f_\lambda(z)=\overline \lambda f(\lambda z)$. It is easy to check that $f_\lambda\in \mathcal K_H^0$ for all $|\lambda| =1$.

\par\smallskip
Rotations are particular cases of holomorphic automorphisms of the unit disk. Other examples are the functions $\varphi_a$, for $a\in\D$, given by
\[
\varphi_a(z)=\frac{a-z}{1-\overline a z}\,,\quad z\in\D\,.
\]

It is not difficult to check that if $f=h+\overline g \in \mathcal K_H^0$  and $a\in\D$, the function
\[
f_a=h_a+\overline{g_a}=\frac{f\circ\varphi_a-f(a)}{h'(a)(1-|a|^2)}\in \mathcal K_H\,.
\]
Since the dilatation $\widehat{\omega_a}=g'_a/h'_a=\mu_a(\omega\circ\varphi_a)$ , where $\omega$ is the dilatation of $f$ and $\mu_a=h'(a)/\overline{h'(a)}$,  does not necessarily fix the origin, one further transformation is needed to get a function in $\mathcal K_H^0$. Concretely, given a function $f=h+\overline g$ in this class with dilatation $\omega$, the new function
\begin{equation}\label{eq-trans}
F_a=\frac{f_a-\overline{\widehat{\omega_a}(0)f_a}}{1-|\omega_a(0)|^2}\in\mathcal K_H^0\,.
\end{equation}
We refer the reader to \cite[p. 79]{Dur-harm} for a proof. A straightforward calculation shows that if we write the function $F_a$ in \eqref{eq-trans} as $F_a=H_a+\overline G_a$, then
\begin{equation}\label{eq-Ha}
H_a(z)=\frac{\left(h\left(\varphi_a(z)\right)-h(a)\right)-\overline{\omega(a)}\left(g\left(\varphi_a(z)\right)-g(a)\right)}{h'(a)(1-|a|^2)(1-|\omega(a)|^2)}
\end{equation}
and 
\[
G_a(z)=\frac{\left(g\left(\varphi_a(z)\right)-h(a)\right)-\omega(a)\left(h\left(\varphi_a(z)\right)-h(a)\right)}{\overline{h'(a)}(1-|a|^2)(1-|\omega(a)|^2)}\,,\quad z\in\D\,.
\]
The dilatation $\omega_a$ of $F_a$ is, for $z$ in the unit disk,
\begin{equation}\label{eq-dilatation}
\omega_a(z)=-\mu_a(\varphi_{\omega(a)}\circ\omega\circ\varphi_a(z))= -\frac{h'(a)}{\overline{h'(a)}}\frac{\omega(a)-\omega(\varphi_a(z))}{1-\overline{\omega(a)}\omega(\varphi_a(z))}\,.
\end{equation}
\section{Proof of Theorem~\ref{thm-main}}\label{sec-main}

\begin{pf}
Let $f=h+\overline g\in\mathcal K_H^0$ have dilatation $\omega$. We start by proving that the inequality 
\[
|h'(z)|\leq \frac{1}{(1-|z|)^3}
\]
holds for $z\in\D$.

\par\smallskip
To do so, let us fix the real values $\alpha$ and $\beta$ provided by Theorem~\ref{thm-CSS2} so that 
\[
{\rm Re\,}\left( e^{i(\alpha+\beta)} \left(h'(z) +e^{-2i\alpha}g'(z)\right)\left(1-e^{-2i\beta}z^2\right)\right)\geq 0\,.
\]
A direct application of Lemma~\ref{lem-q} shows that there exists an analytic function $\delta$ in the unit disk with $\delta(0)=0$ and $\delta(\D)\subset\D$ such that for $z\in\D$,
\[
h'(z)+e^{-2i\alpha}g'(z)=\frac{1+e^{-2i(\alpha+\beta)}\delta(z)}{1-\delta(z)}\frac{1}{1-e^{-2i\beta}z^2}\,. 
\]
As a consequence of the triangle inequality and de Schwarz lemma we have, from the previous identity, that for some $\alpha\in\R$,
\begin{equation}\label{eq-1}
|h'(z)+e^{-2i\alpha}g'(z)|\leq \frac{1+|z|}{1-|z|}\frac{1}{1-|z|^2}=\frac{1}{(1-|z|)^2}\,. 
\end{equation}

\par\smallskip
On the other hand, by Theorem~\ref{thm-CSS1}, the function $h-e^{-2i\alpha}g$ belongs to the class $\mathcal S$. Therefore, by Theorem~\ref{thm-S},
\begin{equation}\label{eq-2}
|h'(z)-e^{-2i\alpha}g'(z)|\leq \frac{1+|z|}{(1-|z|)^3}\,. 
\end{equation}

\par\smallskip
Hence, using again the triangle inequality, jointly with the estimates \eqref{eq-1} and \eqref{eq-2}, we get
\[
|h'(z)| =\frac{1}{2} |h'(z)+e^{-2i\alpha}g'(z)+h'(z)-e^{-2i\alpha}g'(z)| \leq  \frac{1}{(1-|z|)^3}\,.
\]
This proves the upper bound in \eqref{eq-dist-H}. The sharpness of the bound will be deferred to the end of the proof. However, it may be convenient to stress at this point that if for some $z\in\D\setminus\{0\}$ we have $|h'(z)|=1/(1-|z|)^3$, then there must exist $\alpha\in \R$ such that equality in \eqref{eq-2} holds, which gives that for that value of $\alpha$ the function $h-e^{-2i\alpha}g$ is a rotation of the Koebe function \eqref{eq-koebe}, according to Theorem~\ref{thm-S}.

\par\smallskip
To prove that the upper bound in \eqref{eq-growth-H} is satisfied, we argue as usual, by using the bound for $|h'|$ as follows. Let us fix $z=re^{i\theta}\in \D\setminus\{0\}$, so that $|z|=r\in(0,1)$, and observe that 
\[
h(z)=\int_0^r h'(\rho e^{i\theta})\, e^{i\theta}\, d\rho\,,
\]
since $h(0)=0$. Therefore, 
\[
|h(z)|\leq \int_0^r |h'(\rho e^{i\theta})|\, |d\rho|\leq  \int_0^r\frac{1}{(1-|\rho|)^3}\, |d\rho| = \frac{2r-r^2}{2(1-r^2)}\,.
\]
\par\smallskip
The previous estimate proves the upper bound in \eqref{eq-growth-H}. Equality holds for some point $z$ in the unit disk different from zero if equality holds for the bound of $|h'|$ for all points in the segment joining the origin with $z$. Therefore, as above, we conclude that if equality holds for such $z$, then $h-e^{-2i\alpha}g$ is a rotation of the Koebe function \eqref{eq-koebe} for some $\alpha\in\R$.

\par\smallskip
Let us now prove  the lower bound in \eqref{eq-dist-H}: given $f=h+\overline g\in\mathcal K_H^0$ and given $z\in\D$, the inequality
\[
|h'(z)|\geq \frac{1}{(1+|z|)^3}
\] 
holds. 

\par\smallskip

 To do so, let $f=h+\overline g\in\mathcal K_H^0$ and let $a\in\D$. Consider the function $F_a=H_a+\overline G_a$ as in \eqref{eq-trans}. Since $F_a\in\mathcal K_H^0$, we can argue as in \eqref{eq-1} to get a real value $\alpha$ such that for all $z\in\D$,
\[
|H_a'(z)+e^{-2i\alpha}G_a'(z)|\leq \frac{1}{(1-|z|)^2}\,. 
\]
In particular, if we set $z=a$, we have
\begin{equation}\label{eq-3}
|H_a'(a)+e^{-2i\alpha}G_a'(a)|=|H_a'(a)|\cdot |1+e^{-2i\alpha}\omega_a(a)|\leq \frac{1}{(1-|a|)^2}\,.
\end{equation}
Here, $\omega_a$ is the dilatation of $F_a$ given by \eqref{eq-dilatation}, which satisfies $|\omega_a(a)|=|\omega(a)|$, since the dilatation $\omega$ of any function in $\mathcal K_H^0$ fixes the origin. Therefore, an application of the triangle inequality gives, from \eqref{eq-3},
\begin{equation}\label{eq-eq}
|H_a'(a)|(1-|\omega(a)|)\leq |H_a'(a)|\cdot |1+e^{-2i\alpha}\omega_a(a)|\leq \frac{1}{(1-|a|)^2}\,.
\end{equation}

\par\smallskip
Using the formula for $H_a$ in \eqref{eq-Ha}, it is a straightforward calculation to show that, since  $h'(0)=1$, $g'(0)=0$, $\varphi_a(a)=0$, and $|\varphi'(a)|=1/(1-|a|^2)$,
\begin{equation}\label{eq-eq2}
|H_a'(a)|=\frac{1}{|h'(a)| (1-|a|^2)^2 (1-|\omega(a)|^2)}\,.
\end{equation}
Therefore, setting \eqref{eq-eq2} in \eqref{eq-eq}, we conclude that for any $a$ in the unit disk,
\begin{eqnarray*}
\frac{1-|\omega(a)|}{|h'(a)| (1-|a|^2)^2 (1-|\omega(a)|^2)}&=&\frac{1}{|h'(a)| (1-|a|^2)^2 (1+|\omega(a)|)}\\
& \leq& \frac{1}{(1-|a|)^2}\,.
\end{eqnarray*}
This is equivalent to 
\begin{eqnarray}\label{eq-dilat}
\nonumber |h'(a)|&\geq& \frac{(1-|a|)^2}{(1-|a|^2)^2 (1+|\omega(a)|)}\\
&=& \frac{1}{(1+|a|)^2(1+|\omega(a)|))}\geq  \frac{1}{(1+|a|)^3}\,,
\end{eqnarray}
since $\omega$ is an analytic function in the unit disk that satisfies the hypotheses of ht Schwarz lemma and therefore, the inequality $|\omega(a)|\leq |a|$ holds. This proves the lower bound in \eqref{eq-dist-H}.

\par\smallskip
For $0<r<1$, let $m_r=\min_{|z|=r} |h(z)|$. Since $h$ is in the class $\mathcal S$, by Theorem~\ref{thm-CSS1}, we get that $m_r>0$ and by the definition of $m_r$ we have that the Euclidean disk centered at the origin and with radius $m_r$ is contained in the set $S=\{h(z)\colon |z|\leq r\}$. Moreover, by the continuity of $h$ in $\D$, there must exists a point $z_0$ of modulus $r$ such that $|h(z_0)|=m_r$, so that the segment $s$ which joints $0$ and $h(z_0)$ is contained in $S$. Let $C$ be the pre-image under $h$ of the segment $s$. Using again that $h$ is univalent in the unit disk, we have that $C$ is an simple analytic curve  which joins $0$ and $h(z_0)$. We then have, since $h'(\z)\, d\z$ has constant argument for $\z\in C$, that using that the lower bound in \eqref{eq-dist-H} holds,
\begin{eqnarray}\label{eq-final}
|h(z_0)|&=&\left|\int_C h'(\z))\, d\z\right|= \int_C |h'(\z))|\, |d\z|\geq \int_0^r \frac{1}{(1+\rho)^3}\, d\rho \\
\nonumber &=& \frac{2r+r^2}{2(1+r)^2}\,. 
\end{eqnarray}
This shows that if $|z|=r$, then $$|h(z)|\geq |h(z_0)|\geq \frac{2r+r^2}{2(1+r)^2}\,,$$ which proves the lower bound in \eqref{eq-growth-H}. 

\subsection*{Sharpness}
To finish the proof of Theorem~\ref{thm-main}, we are to show that equality in any of the inequalities in either \eqref{eq-growth-H} or \eqref{eq-dist-H} holds for some $z\in\D\setminus\{0\}$ for a function $f\in\mathcal K_H^0$, then $f$ is a rotation of the function $L$ defined in \eqref{eq-L}. To do so, we will use the following proposition. We have not been able to find any explicit reference of this result, so that we include a proof for the sake of completeness.  Notice that the function $L$ in \eqref{eq-L} satisfies the condition $H-G=k$.

\begin{thm}\label{thm-Koebe}
Let $f=h+\overline g\in\mathcal K_H^0$. If there exists a real number $\mu$ such that $h-e^{-2i\mu} g$ equals a rotation of the Koebe function $k$ in \eqref{eq-koebe}, then $f$ is a rotation of the half-plane harmonic mapping $L$.
\end{thm}
\begin{pf}
Let $\lambda=e^{it}$, $t\in\R$, be such that $h(z)-e^{-2i\mu} g(z)=e^{-it} k(e^{it}z)$. The rotation $f_\mu=h_\mu+\overline{g_\mu}$ of $f$, defined by $f_\mu(z)=e^{i\mu}f(e^{-i\mu}z)$ satisfies 
\[
h_\mu'(z)-g_\mu'(z)=h'(e^{-i\mu}z)- e^{-2i\mu}g'(e^{-i\mu}z)=k'(e^{i(t-\mu)} z)\,.
\]
This implies that $h_\mu-g_\mu$ equals a rotation of the Koebe function. Concretely,  $ h_\mu(z)-g_\mu(z)=e^{-i(\mu-t)}k(e^{i(t-\mu)} z)$, $z\in\D$. Since $f_\mu$ is a rotation of a function in $\mathcal K_H^0$, it belongs to $\mathcal K_H^0$ as well. Hence, according to Theorem~\ref{thm-CSS1}, $h_\mu-g_\mu$ must be convex in the $0$ direction. 

\par\smallskip 
The Koebe function $k$ maps the unit disk onto the complex plane minus the slit along the negative real axes from $-\infty$ to $-1/4$. Therefore, $k$ is convex only in the horizontal direction (that is, only in the $0$ direction) and hence, the rotation $k_\lambda$ of the function $k$ for $\lambda=e^{i\theta}$, given by $k_\lambda(z)=e^{i\theta}k(e^{-i\theta}z)$, is convex only in the $\theta$ direction for $\theta\in[0,\pi)$, and only in the $2\pi-\theta$ direction for $\theta\in[\pi, 2\pi)$. This gives only two possibilities: either $h_\mu-g_\mu=k$ or $h_\mu-g_\mu=\widehat k$, where $\widehat k(z)=-k(-z)$. The latter case can be reduced to the former by applying a new rotation to the function $f_\mu$. Concretely, $\widehat{f_\mu}=-f_\mu(-z)$. Therefore, we may assume, without loss of generality, that if a given function satisfies the hypotheses in the theorem, then some of its rotations that we again denote by $f=h+\overline g$ in order not to burden the notation, satisfies $h-g=k$. The proof of the proposition will be completed if we check that any such function must be the function $L$.
\par\smallskip

To check that any convex harmonic mapping $f=h+\overline g$ in $\mathcal K_H^0$ with $h-g=k$ is the harmonic half-plane mapping, we proceed as follows. 
\par\smallskip
Recall that $h-g\in \mathcal S$. Hence, by the Bieberbach Theorem~\ref{thm-Bie}, $$|h''(0)-g''(0)|\leq 4\,.$$ On the other hand, we have, by Theorem~\ref{thm-CSS3}, that  $|h''(0)|\leq 3$ and $|g''(0)|\leq 1$. Therefore, we obtain
\[
4=k''(0)={\rm Re\,} k''(0)= {\rm Re\,}\left(h''(0)-g''(0)\right)\leq |h''(0)|+|g''(0)|\leq 4\,.
\]
That is, in particular, ${\rm Re\,} g''(0)=-1$. 

\par\smallskip
The dilatation $\omega$ of $f$ is given by $\omega(z)=g'(z)/h'(z)$ for $z$ in the unit disk. The conditions $\omega(0)=0$ and $\omega(\D)\subset \D$ are satisfied. Therefore, by the Schwarz lemma, $|\omega'(0)|\leq 1$, with equality if and only if $\omega(z)=\lambda z$ for some $\lambda=1$. It is easy to check that $\omega'(0)=g''(0)$. Hence, $\omega(z)=-z$.

\par\smallskip
The solution of the system of equations $h-g=k$ and $g'(z)=-zh'(z)$, $h(0)=g(0)=0$, is the function $L$. This ends the proof.
\end{pf}
 
 \par\smallskip
 Let us now finish the proof of Theorem~\ref{thm-main}. As pointed out above, if a given function $f=h+\overline g\in\mathcal K_H^0$ satisfies that $$|h(z)|=\frac{2|z|-|z|^2}{2(1-|z|)^2}$$ for some $z\in\D\setminus\{0\}$, then, necessarily, $|h'(z)|=1/(1-|z|)^3$. In order that this latter condition to hold, the inequality in \eqref{eq-2} must be satisfied at $z$. According to Theorem~\ref{thm-S}, $h-e^{i\theta}g$ is a rotation of the Koebe function. Therefore, by Theorem~\ref{thm-Koebe}, $f$ is a rotation of the half-plane mapping $L$ and hence, $h$ is a rotation of the function $H$ as in \eqref{eq-H}.

\par\smallskip
Regarding the sharpness of the lower bounds in \eqref{eq-growth-H} and \eqref{eq-dist-H}, according to \eqref{eq-final}, we have that if the lower bound in  \eqref{eq-growth-H} is attained by some $a\in\D\setminus\{0\}$, then so is the lower bound in  \eqref{eq-dist-H} for the same $a$. Hence, it only remains to analyze the case of equality in the lower bound of  \eqref{eq-dist-H}. 

\par\smallskip
Assume then that for some $a\in\D\setminus\{0\}$, the analytic part $h$ of a function $f=h+\overline g\in\mathcal K_H^0$ satisfies
\[
|h'(a)|=\frac{1}{(1+|a|)^3}\,.
\]

In this case, equality must hold in \eqref{eq-dilat}. So that $\omega$ is a rotation and we can write $\omega(z)=\lambda z$ for some $|\lambda|=1$. 

\par\smallskip
In this case, the analytic function $h-\nu g$, where $\nu=\overline \lambda\overline a /|a|$ belongs to the class $\mathcal S$ and  satisfies that
\[
|h'(a)-\nu g'(a)|= |h'(a)||1-\nu \omega(a)|=\frac{1-|a|}{(1+|a|)^3}\,.
\]
Therefore, by Theorem~\ref{thm-S}, $h-\nu g$ is a rotation of the Koebe function. By Theorem~\ref{thm-Koebe}, $f$ is a rotation of the half-plane mapping $L$ and hence, $h$ is a rotation of the function $H$ as in \eqref{eq-H}. This completes the proof of Theorem~\ref{thm-main}.
\end{pf}

\subsection*{Further remarks}
A careful reader may have noticed some direct consequences of both Theorem~\ref{thm-main} and Theorem~\ref{thm-Koebe} that we now list.

\par\smallskip
The first one is that not any function in the class $\mathcal S$ is allowed to be the analytic part of a function $f=h+\overline g\in\mathcal K_H^0$, since the bounds for the growth and distortion for the function $h$ given by Theorem~   \ref{thm-main} do not reach the values of the bounds for functions in the class $\mathcal S$ provided by Theorem~\ref{thm-S}. In particular, the Koebe function \eqref{eq-koebe} cannot be the analytic part of a function $f\in \mathcal K_H^0$, a result that also follows from Theorem~\ref{thm-CSS3}.

\par\smallskip
More generally, the arguments used to prove Theorem~\ref{thm-Koebe} show that, indeed, if for some $a$ in the unit disk a function $f=h+\overline g \in\mathcal K_H^0$ satisfies that $h+ag$ is a rotation of the Koebe function, then $|a|=1$ and $f$ is a rotation of the harmonic half-plane mapping. 

\section*{Acknowledgements}
It is a pleasure to thank Iason Efraimidis and Rodrigo Hern\'andez for their careful reading of the first handwritten draft of the proofs of the  theorems in this paper and for their comments that have helped me to improve the exposition.

\end{document}